\documentclass[12pt,letterpaper,reqno]{amsart}

\usepackage{amssymb,latexsym, amsmath, amsxtra, amssymb,amscd,
verbatim} \usepackage{color}
\usepackage[all]{xy}
\usepackage[dvips]{graphics}

\theoremstyle{plain}
\newtheorem{theorem}{Theorem}[section]

\newtheorem{lemma}[theorem]{Lemma}
\newtheorem{proposition}[theorem]{Proposition}
\theoremstyle{definition}

\theoremstyle{remark}

\makeatletter
\def\Ddots{\mathinner{\mkern1mu\raise\p@
\vbox{\kern7\p@\hbox{.}}\mkern2mu
\raise4\p@\hbox{.}\mkern2mu\raise7\p@\hbox{.}\mkern1mu}}
\makeatother

\numberwithin{equation}{section}
\numberwithin{theorem}{section}
\numberwithin{table}{section}
\numberwithin{figure}{section}

\newcommand{\C}{\mathbb C}

\newcommand{\R}{\mathbb R}
\newcommand{\Z}{\mathbb Z}

\renewcommand{\H}{\mathcal H}

\newcommand{\cj}[1]{\overline{#1}}
\newcommand{\abs}[1]{\left| #1 \right|}
\newcommand{\cbr}[1]{\left\{ #1 \right\}}

\def\({\left(}
\def\){\right)}

\newcommand{\A}{{\mathbb A}}

\newcommand{\GL}{{\rm GL}}

\newcommand{\GSp}{{\rm GSp}}

\newcommand{\PGSp}{{\rm PGSp}}

\newcommand{\SSp}{{\rm Sp}}

\newcommand{\mat}[4]{{\setlength{\arraycolsep}{0.5mm}\left[
\begin{array}{cc}#1&#2\\#3&#4\end{array}\right]}}

\begin{document}
\title[Testing a functional equation]{Testing the functional equations of a high-degree Euler product.}
\author[Farmer, Ryan, and Schmidt]{David W. Farmer, Nathan C. Ryan, and Ralf Schmidt}
\date{\today}


\thispagestyle{empty}
\vspace{.5cm}
\begin{abstract}

We study the L-functions associated to Siegel modular forms (equivalently, automorphic representations of ${\rm GSp}(4,\mathbb{A}_{\mathbb{Q}})$) both theoretically and numerically.  For the L-functions
of degrees 10, 14, and 16 we perform representation theoretic
calculations to cast the Langlands L-function in classical terms.
We develop a precise notion of what it means to test a
conjectured functional equation for an L-function, and we
apply this to the degree 10 adjoint L-function associated to a Siegel
modular form.

\end{abstract}

\address{
{\parskip 0pt
American Institute of Mathematics\endgraf
farmer@aimath.org\endgraf
\null
Department of Mathematics, Bucknell University\endgraf
nathan.ryan@bucknell.edu \endgraf
\null
Department of Mathematics, University of Oklahoma \endgraf
rschmidt@math.ou.edu
}
  }

\maketitle

\section{Introduction}

L-functions are special functions that arise in representation theory
and in several areas of number theory.
From the viewpoint of analytic number theory, L-functions are Dirichlet series
with a functional equation and an Euler product.  From the point of view of representation theory,
L-functions arise from automorphic representations of
a reductive group over the adeles of a number field.

The two points of view offer distinct benefits.
Representation theory, via the Langlands program~\cite{Gelbart},
provides a framework for understanding how L-functions arise,
as well as the connections between various mathematical objects.
L-functions considered as objects of analytic number theory are suitable for concrete
exploration and testing of conjectures, for example they
can be evaluated on a computer to check the Riemann hypothesis.
Unfortunately, it can be quite difficult to translate
Langlands L-functions into this setting, which limits the
ability to do explicit calculations and test conjectures.  In this paper
we make such a translation and perform computer calculations with the
results: for a particular Siegel modular form $F$, we calculate factors
$L_p(s,F,\rho)$ and $\varepsilon_p(s,F,\rho)$ ($p\leq \infty)$ for
six choices of $\rho$ (dimensions 4, 5, 10,14, 16) using the Langlands
parameterization of the discrete series.  Using these calculations we
provide numerical evidence that the L-function of degree 10 satisfies
a functional equation.


The L-functions we consider here arise (in the classical setting)
from holomorphic Siegel modular forms on $\SSp(4,\Z)$, see Section~\ref{sec:Siegel}.
The same L-functions arise 
from automorphic representations of
$\PGSp(4,\mathbb{A})$, see Section~\ref{sec:Langlands} (as in \cite{AS1}, at the archimedean place such an automorphic representation is a holomorphic discrete series representation and at the nonarchimedean places it is a spherical principal series representation).
The Langlands program predicts the existence of an infinite list of L-functions associated to a Siegel modular form.  In our particular case the first two L-functions
are known as the spinor and the standard L-function, and have
degree~4 and~5, respectively.  Due to
Andrianov \cite{An}, Shimura \cite{Shimura}, B\"ocherer \cite{Boecherer}
and others,
these L-functions are fairly well-understood;
for Siegel modular forms on $\SSp(4,\Z)$, they are known to be entire functions that satisfy a functional
equation.  

The next case is the adjoint L-function, which has degree~10.
It has not been shown that this L-function is entire and satisfies
a functional equation. 
(The automorphic representations associated to
holomorphic Siegel modular forms are not generic, so
a technique as in~\cite{Ginzburg} is not applicable.)
Providing evidence for the conjectured functional equation,
via a computer calculation, is one of the goals of
this paper.
This is made precise in Theorem~\ref{thm:testFE},
which gives a test for the functional equation and quantifies
the probability that the test will yield a false positive.


A substantial part of this paper is a translation from
the perspective of representation theory to the viewpoint
of analytic number theory.  Selberg~\cite{Sel} gave a set of axioms for what is
now called the ``Selberg class'' of L-functions.
We will call L-functions in this class ``Selberg L-functions,'' and are to be compared with Langlands L-functions -- those that arise from automorphic representations.  It is a standard conjecture that all Selberg L-functions are Langlands L-functions and that all primitive Langlands L-functions are Selberg L-functions.  In this paper we translate Langlands L-functions into Selberg L-functions.  We now describe Selberg L-functions in more detail.

A Selberg L-function $L(s)$ is given by a Dirichlet series $L(s) = \sum_{n=1}^\infty \frac{a_n}{n^s}$ where $a_1=1$ and the series converges in some half-plane.  We assume a \emph{Ramanujan bound} on the coefficients: $a_n = O(n^\varepsilon)$ for
any $\varepsilon>0$.
Moreover, it has a meromorphic continuation to the
whole complex plane with at most finitely many poles,
all of which are on the line~$\Re(s)=1$.  $L(s)$ can be written as an \emph{Euler product} $L(s)=\prod_p L_p(p^{-s})^{-1}$
where the product is over the primes, and $L_p$ is a polynomial
with $L_p(0)=1$.  Additionally, there exist $Q>0$,
positive real numbers $\kappa_1$,\ldots,$\kappa_n$,
and complex numbers with non-negative real part
$\mu_1$,\ldots,$\mu_n$, 
such that $\Lambda(s):= Q^s \prod \Gamma(\kappa_j s + \mu_j) \cdot L(s)$
is meromorphic with poles only arising from the poles of $L(s)$ and
satisfies the functional equation
$\Lambda(s) = \varepsilon \overline{\Lambda(1-\overline{s})}$ where
$|\varepsilon|=1$.  The number $d=2\sum \kappa_j$ is the \emph{degree}
of the L-function.  An alternate way of thinking about the Ramanujan bound is that $|a_p|\leq d$.  

 Indeed, for the L-functions considered here there do
not exist results in the literature which would allow a non-expert
to translate the L-function data from representation theoretic
language into a form involving a Dirichlet series.
Thus,
we give a brief introduction to the aspects of the Langlands
program that are relevant to our calculations, and 
describe how to translate
the Langlands L-functions we consider here into Selberg L-functions.  
The results of those calculations are summarized in
Proposition~\ref{prop:45101416}.


In Section~\ref{sec:Siegel}
we describe the Siegel modular forms which give rise to the
L-functions considered here, and we describe these L-functions
in the classical language.  In Section~\ref{sec:Langlands} we describe Langlands L-functions
and how to translate the degree 10 L-function considered here into a
form which can
be evaluated on a computer.  In Section~\ref{sec:checkFE}
we provide evidence for the conjectured functional equations and also 
briefly address the problem of accurately evaluating L-functions for which
only a few of the local factors in the Euler product are known.  We also provide criteria to measure the strength of the evidence.

We would like to thank Sharon Garthwaite for her helpful comments on the paper and the referee for helping improve the paper.

\section{Siegel modular forms and their L-functions}\label{sec:Siegel}
We recall the definition and main properties of Siegel
modular forms on $\SSp(2n,\Z)$ and we describe the two simplest
L-functions associated to them.

%
%
%
\subsection{Siegel modular forms}\label{ssec:siegel}
Let $0_n$ be the zero matrix, $E_n:=\left(\begin{smallmatrix}& & 1\\
& \Ddots &\\
1 & & \\
\end{smallmatrix}\right)$ and $J_n:=\left(\begin{smallmatrix} 0_n& E_n\\ -E_n&0_n\end{smallmatrix}\right)$.  Denote the group of symplectic similitudes by $
\GSp^+(2n,\R):=\{\alpha \in \text{GL}(2n,\R) : {}^t\alpha J_n\alpha = r(\alpha)J_n, r(\alpha) \in \R, r(\alpha)>0)\}$,
where $r(\alpha)$ is called the similitude of $\alpha$.
Define the Siegel modular group of genus $n$ by $
\SSp(2n,\Z) := \{\gamma \in \GSp^+(2n,\R) \cap M(2n,\Z) : r(\gamma)=1\}$.
Let $\H^n: = \{Z = X+iY: X,Y \in M(n,\R), {}^tZ=Z, Y>0\}$ denote the Siegel upper half space, that is, symmetric matrices in $M(n,\C)$ with positive definite imaginary part.

Recall, a holomorphic function $F: \H^n \to \C$ is a Siegel modular form
of genus $n$ and weight $k$ 
if for all $\alpha=\left(\begin{smallmatrix}A &B \\C &
D\end{smallmatrix}\right) \in \SSp(2n,\Z)$ it satisfies the
transformation property 
\[ F(Z)= (F|_k\alpha)(Z) :=
r(\alpha)^{nk-n(n+1)/2}\det(CZ+D)^{-k}F((AZ+B)(CZ+D)^{-1}). 
\]
If $n=1$ then $F$ must satisfy an additional
growth condition.  

We shall denote the space of weight $k$
genus $2$ Siegel modular forms by $M_k(\SSp(4,\Z))$.

In genus 2, we can express the expansion of a Siegel cusp form as
\[
F(Z)=\sum_{\substack{r,n,m\in\mathbf{Z}\\r^2-4mn< 0\\n,m\geq 0}} a_F(n,r,m) q^n\zeta^r q'^n
\]
where $[n,r,m]$ is the positive definite binary quadratic form $nX^2+rXY+mY^2$ of discriminant $r^2-4mn$ and $q=e^{2\pi i z}~(z\in\H^1)$, $q' = e^{2\pi i \omega}~(\omega\in\H^1)$, and $\zeta = e^{2\pi i \tau}~(\tau\in\C)$.  In particular, we are examining L-functions associated to modular forms not in the Maass space, i.e., not  in the image of the Saito Kurokawa lift, i.e., whose L-functions are primitive (at least conjecturally).  
The first such form occurs in weight $k=20$ and is computed in \cite{Skoruppa}.

There is a theory of Hecke operators acting on the
space of Siegel modular forms; we denote the $n$th Hecke operator by $T(n)$.  The Hecke eigenvalues for a Siegel modular form can be computed explicitly from its Fourier coefficients, but this is computationally expensive.  Let $F$ be a Hecke eigenform; i.e., suppose that for each $n$ there exists a $\lambda_F(n)$ so that $F|T(n) = \lambda_F(n) F$.  For the weight
20 Siegel cusp form $F$ that is not a Saito-Kurokawa lift, we will use $\lambda_F(p)$ and $\lambda_F(p^2)$ for $p\leq 79$ in our experiments.  These data are computed in~\cite{KohnenKuss}.


\subsection{L-functions associated to Siegel modular forms}

There are two well-known L-functions attached to Siegel
modular forms on $\SSp(4,\Z)$, called the spinor L-function
and the standard L-function.  These have been studied by
Andrianov \cite{An}, Shimura \cite{Shimura}, B\"ocherer \cite{Boecherer}
and others.  Formulas for those L-functions in the genus~2 case
are given in Proposition~\ref{prop:45101416}.  

To each genus $2$ eigenform, $F$, one can associate, for each prime $p$, a triple $(\alpha_{0,p},\alpha_{1,p},\alpha_{2,p})$ of nonzero complex numbers -- it is in these terms that our L-functions are expressed.  The entries of the triple are called the Satake parameters of $F$.  

In genus 2 it is rather straightforward to compute the Satake parameters of a form, given the Hecke eigenvalues of the form.  By using an explicit description of Satake isomorphism as found, for example, in \cite{RyanShemanske}, write the Euler factor of the spinor L-function as a polynomial whose coefficients are in terms of $\lambda_F(p)$ and $\lambda_F(p^2)$.  To compute the Satake parameters, one finds the roots of this polynomial.

We rescale the Satake parameters to have the
normalization $|\alpha_j|=1$,  $ \alpha_0^2\alpha_1\alpha_2 =1$,
which is possible since the Ramanujan bound is a theorem for Siegel modular forms \cite{Weissauer}.
This
corresponds to a simple change of variables in the L-functions, so that all our L-functions satisfy
a functional equation in the standard form~$s\leftrightarrow 1-s$.


%

In Section~\ref{sec:Langlands} we describe the procedure for determining the L-functions associated to an automorphic
representation, and give reasonably complete details for the
spinor, standard, and adjoint L-functions of genus~2 Siegel modular
forms.  The results of those calculations are summarized in the
following proposition.  For completeness we also include the results of similar computations carried out for two more L-functions: ones associated to a specific 14-dimensional representation and 16-dimensional representation.

\begin{proposition}\label{prop:45101416}  Suppose $F\in M_k(\SSp(4,\Z))$ be a Hecke eigenform.  Let $\alpha_{0,p}$, $\alpha_{1,p}$,
$\alpha_{2,p}$ be the Satake parameters of $F$ for the prime~$p$,
where we suppress the dependence on $p$ in the formulas below.  For $\rho\in\{\mathrm{spin},\mathrm{stan},\mathrm{adj},\rho_{14},\rho_{16}\}$ we have the L-functions $L(s,F,\rho):=\prod_{p \text{ prime}} Q_p(p^{-s},F,\rho)^{-1}$ where

\begin{align}\label{eqn:EPs}
Q_p(X,F,\mathrm{spin}):=\mathstrut &
{(1-\alpha_0 X)(1-\alpha_0\alpha_1 X)
 (1-\alpha_0\alpha_2 X)(1-\alpha_0\alpha_1\alpha_2 X)},\nonumber\\
 Q_p(X,F,\mathrm{stan}):=\mathstrut &(1-X)(1-\alpha_{1}X)(1-\alpha_{1}^{-1}X)
(1-\alpha_{2}X)(1-\alpha_{2}^{-1}X),\nonumber\\
Q_p(X,F,\mathrm{adj}):=\mathstrut &(1- X)^2(1-\alpha_1 X)(1-\alpha_1^{-1} X)
 (1-\alpha_2 X)(1-\alpha_2^{-1} X)\nonumber\\
 &\qquad(1-\alpha_1\alpha_2 X)(1-\alpha_1^{-1}\alpha_2 X)
  (1-\alpha_1\alpha_2^{-1} X)(1-\alpha_1^{-1}\alpha_2^{-1} X),\nonumber\\
  Q_p(X,F,\rho_{14}):=\mathstrut&(1- X)^2(1-\alpha_1 X)(1-\alpha_2 X)(1-\alpha_1^{-1} X)(1-\alpha_2^{-1} X)\\
 &\qquad(1-\alpha_1^2 X)(1-\alpha_2^2 X) (1-\alpha_1^{-2} X)(1-\alpha_2^{-2} X)\nonumber\\
 &\qquad(1-\alpha_1\alpha_2 X)(1-\alpha_1\alpha_2^{-1} X) (1-\alpha_1^{-1}\alpha_2 X)(1-\alpha_1^{-1}\alpha_2^{-1} X),\mathrm{ and}\nonumber\\
  Q_p(X,F,\rho_{16}):=\mathstrut
&(1-\alpha_0 X)^2(1-\alpha_0\alpha_1 X)^2
  (1-\alpha_0\alpha_2 X)^2(1-\alpha_0\alpha_1\alpha_2 X)^2\nonumber\\
 &(1-\alpha_0\alpha_1^{-1} X)(1-\alpha_0\alpha_2^{-1} X)
  (1-\alpha_0\alpha_1^2 X)(1-\alpha_0\alpha_2^2 X)\nonumber\\
 &(1-\alpha_0\alpha_1\alpha_2^{-1} X)(1-\alpha_0\alpha_1^{-1}\alpha_2 X)
  (1-\alpha_0\alpha_1^2\alpha_2 X)(1-\alpha_0\alpha_1\alpha_2^2 X).\nonumber
 \end{align}
give the L-series of, respectively, the spinor, standard, adjoint, degree 14 and degree 16 L-functions.  These L-functions satisfy the functional equations (the last three conjecturally satisfy the functional equations):
\begin{align}\label{eqn:FEs}
\Lambda(s,F,\mathrm{spin}):=\mathstrut &
     (4\pi^2)^{-s} \Gamma(s+\tfrac12)\Gamma(s+k-\tfrac32)L(s,F,\mathrm{spin})\nonumber \\
=\mathstrut & (-1)^k \Lambda(1-s,F,\mathrm{spin}),\nonumber\\
\Lambda(s,F,\mathrm{stan}):=\mathstrut &  (4 \pi^{5/2})^{-s} \Gamma(\tfrac12 s)\Gamma(s+k-2)\Gamma(s+k-1) L(s,F,\mathrm{stan}) \nonumber\\
=\mathstrut & \Lambda(1-s,F,\mathrm{stan})\nonumber,\\
\Lambda(s,F,\text{adj}):=\mathstrut &  (16 \pi^{5})^{-s}
\Gamma(\tfrac12 (s+1))^2\Gamma(s+1)\cr & \times\Gamma(s+k-2)\Gamma(s+k-1)
\Gamma(s+2k-3)  L(s,F,\mathrm{adj}) \cr =\mathstrut  &\Lambda(1-s,F,\mathrm{adj}),\\
\Lambda(s,F,\rho_{14}):=\mathstrut &  
(2^6 \pi^7)^{-s}\Gamma(s/2)^2
\Gamma(s+1) 
\Gamma(s+k-2)
\Gamma(s+k-1)
\nonumber\\
 &\qquad\times
\Gamma(s+2k-4)
\Gamma(s+2k-3)
\Gamma(s+2k-2)
L(s,F,\rho_{14}) \cr
=\mathstrut & \Lambda(1-s,F,\rho_{14}),\mathrm{ and}\nonumber\\
 \Lambda(s,F,\rho_{16}):=\mathstrut &
(2^8\pi^8)^{-s}
\Gamma\big(s+\tfrac12\big)^2 
\Gamma\big(s+k-\tfrac52\big)
  \Gamma\big(s+k-\tfrac32\big)^2\,
  \Gamma\big(s+k-\tfrac12\big)\cr
 &\qquad \times
 \Gamma\big(s+2k-\tfrac52\big)\Gamma\big(s+2k-\tfrac72\big)
L(s,F,\rho_{16}) \cr
=\mathstrut & -\Lambda(1-s,F,\rho_{16}).\nonumber
\end{align}
\end{proposition}

\section{Langlands L-functions}\label{sec:Langlands}

The Euler products in the previous section arise as
Langlands $L$-functions attached to automorphic representations of $\GSp(4,\A)$
generated by the Siegel modular form $F$, see \cite{AS1}. In general, this procedure
involves the local Langlands correspondence, which is now a theorem for $\GSp(4)$;
see \cite{Gan}. However, since we are only interested in full level Siegel modular
forms, the mechanism simplifies considerably. We shall briefly describe how to obtain
the local factors in the non-archimedean and the archimedean case.
\subsubsection*{The non-archimedean factors}
Let $\alpha_0,\alpha_1,\alpha_2$ be the Satake parameters of $F$ at a finite place $p$,
normalized as above, so that $\alpha_0^2\alpha_1\alpha_2=1$. These determine
a semisimple conjugacy class in the dual group $\SSp(4,\C)$, represented by the
diagonal matrix
\begin{equation}\label{GSP4Apeq}
 A_{\pi_p}=
\rm{Diag}\left(\alpha_0, \alpha_0\alpha_1,
\alpha_0\alpha_2, \alpha_0\alpha_1\alpha_2\right).
\end{equation}
(One has to carefully go through the definitions of the local Langlands correspondence
to see this; see Sect.\ 2.3 and 2.4 of \cite{RS}.) Let
$\rho:\SSp(4,\C)\rightarrow\GL(n,\C)$ be a finite-dimensional
representation of the dual group. The local $L$-factor attached to
the data $\alpha_0,\alpha_1,\alpha_2$ and $\rho$ is given by
\begin{equation}\label{nonarchLfaceq}
 L_p(s,F,\rho)=\frac1{\det(1-p^{-s}\rho(A_{\pi_p}))}.
\end{equation}
The three smallest non-trivial irreducible representations of $\SSp(4,\C)$ are the
four-dimensional ``spin'' representation (which is simply the inclusion of
$\SSp(4,\C)$ into $\GL(4,\C)$), the five-dimensional ``standard'' representation
(described explicitly in appendix A.7 of \cite{RS}), and the ten-dimensional
adjoint representation ${\rm adj}$ on the Lie algebra $\mathfrak{sp}(4,\C)$.
Calculations show that the resulting $L$-factors are given as follows,
\begin{align*}
 L_p(s,F,{\rm spin})&=Q_p(p^{-s},F,{\rm spin}),\\
 L_p(s,F,{\rm stan})&=Q_p(p^{-s},F,{\rm stan}),\\
 L_p(s,F,{\rm adj})&=Q_p(p^{-s},F,{\rm adj}),
\end{align*}
with the factors $Q_p$ as in the previous section.
There are also corresponding local $\varepsilon$-factors $\varepsilon_p(s,F,\rho)$,
which for unramified representations are all constantly $1$.
\subsubsection*{The archimedean factors}
The real Weil group $W_\R$ is given by $W_\R=\C^\times\sqcup j\C^\times$ with the rules
$j^2=-1$ and $jcj^{-1}=\bar c$ (see \cite{Ta} (1.4.3)). The
commutator subgroup is $S^1\subset\C^\times$, the set of complex
numbers with absolute value $1$. There is a reciprocity law isomorphism
\begin{align}\label{realreciprocityeq}
 r_\R:\:\R^\times&\stackrel{\sim}{\longrightarrow}W_\R^{\rm ab},\\
  -1&\longmapsto jS^1,\nonumber\\
  \R_{>0}\ni x&\longmapsto\sqrt{x}S^1.\nonumber
\end{align}
Let $|\cdot|$ be the usual absolute value on $\R$, and let $\|\cdot\|$ be the
character of $W_\R$ defined by the commutativity of the
following diagram,
\begin{equation}\label{omegadiagdefeq}
 \xymatrix{
 \R^\times\ar[r]^\sim\ar[rd]_{|\,|}&W_\R^{\rm ab}&W_\R\ar[l]\ar[ld]^{\|\cdot\|}\\
 &\C^\times}
\end{equation}
(see \cite{Ta} (1.4.5)). Hence, $\|z\|=|z|^2$ for $z\in\C^\times$,
where $|\cdot|$ denotes the usual absolute value on $\C$. The
character $\|\cdot\|^s$ is denoted by $\omega_s$, for a complex
number $s$ (see \cite{Ta} (2.2)). There are $L$- and $\varepsilon$-factors attached
to characters of $\R^\times$ (see \cite{Ta} (3.1)). The
correspondence between characters of $W_\R$ and characters of
$\R^\times$, and the associated $L$- and $\varepsilon$-factors, are
given in the following table.
\begin{equation}\label{archLLCGL1table}\renewcommand{\arraystretch}{1.4}
 \begin{array}{|c|c|c|c|}
 \hline\text{char.\ of }W_\R&\text{char.\ of }\R^\times
 &\text{L-factor}&\text{$\varepsilon$-factor}\\\hline\hline
 \varphi_{+,t}: z \mapsto |z|^{2t},\;j \mapsto 1&x\mapsto|x|^t
  &\pi^{-(s+t)/2} \Gamma(\frac{s+t}{2})&1\\\hline
 \varphi_{-,t}: z \mapsto |z|^{2t},\;j \mapsto -1&x\mapsto{\rm
sgn}(x)|x|^t   &\pi^{-(s+t+1)/2}
\Gamma(\frac{s+t+1}{2})&i\\\hline  \end{array} \end{equation}
Besides one-dimensional representations, the only other
irreducible representations of $W_\R$ are two-dimensional and indexed
by pairs $(\ell,t)$, where $\ell$ is a positive integer and $t\in\C$. The representation
attached to $(\ell,t)$ is $\varphi_{\ell,t}$, given by
\begin{equation}\label{WR2dimrepeq}
 \varphi_{\ell,t}: re^{i\theta} \longmapsto\mat{r^{2t}
 e^{i\ell\theta}}{}{}{r^{2t} e^{-i\ell\theta}},\quad
   j \longmapsto\mat{}{(-1)^\ell}{1}{}.
\end{equation}
The associated $L$- and $\varepsilon$-factors are
\begin{equation}\label{WR2dimrepLeq}
 L(s,\varphi_{\ell,t})=2(2\pi)^{-(s+t+\ell/2)}\Gamma\Big(s+t+\frac{\ell}{2}\Big),\qquad
 \varepsilon(s,\varphi_{\ell,t})=i^{\ell+1}.
\end{equation}
Now, to a Siegel modular form $F$ of weight $k$ there is attached the four-dimensional
representation of $W_\R$ given by
\begin{equation}\label{weightkparametereq}
 \varphi_{(k-1,k-2)}:=\varphi_{1,0}\oplus\varphi_{2k-3,0}
\end{equation}
(this is really the parameter of a holomorphic discrete series representation with
Harish-Chandra parameter $(k-1,k-2)$, see \cite{Borel}). The image of this
parameter can be conjugated into the dual group $\SSp(4,\C)$. Given a finite-dimensional
representation $\rho:\SSp(4,\C)\rightarrow\GL(n,\C)$, we compose $\rho$ with the
representation (\ref{weightkparametereq}) and obtain an $n$-dimensional
representation of $W_\R$. By \cite{Kn}, this representation can be decomposed into one- and
two-dimensional irreducibles. The product of the $L$-factors (resp.\ $\varepsilon$-factors)
attached to these irreducibles is by definition $L_\infty(s,F,\rho)$
(resp.\ $\varepsilon_\infty(s,F,\rho)$). Calculations show that
\begin{align*}
 L_\infty(s,F,{\rm spin})&=\textstyle4(2\pi)^{-2s-k+1}
  \Gamma(s+\frac12)\Gamma(s+k-\frac32),\\
 L_\infty(s,F,{\rm stan})&=\textstyle2^{-2s-2k+5}\,\pi^{-5s/2-2k+3}\,
  \Gamma(\frac12s)\Gamma(s+k-1)\Gamma(s+k-2),\\
 L_\infty(s,F,{\rm adj})&=\textstyle2^{-4s-3k+9}\pi^{-5s-3k+4}\Gamma(\frac12(s+1))^2
  \Gamma(s+1)\\
  &\qquad\times\Gamma(s+k-2)\Gamma(s+k-1)\Gamma(s+2k-3),
\end{align*}
and
$$
 \varepsilon_\infty(s,F,{\rm spin})=(-1)^k,\qquad
 \varepsilon_\infty(s,F,{\rm stan})=1,\qquad
 \varepsilon_\infty(s,F,{\rm adj})=1.
$$
We see that, up to an irrelevant constant, the archimedean $L$-factors coincide
with the $\Gamma$-factors in \eqref{eqn:FEs}.
\subsubsection*{The global $L$-function}
Having defined all local factors, the global $L$-function attached to $F$ and a
finite-dimensional representation $\rho:\SSp(4,\C)\rightarrow\GL(n,\C)$
is given by
$$
 \Lambda(s,F,\rho)=\prod_{p\leq\infty}L_p(s,F,\rho).
$$
Up to a constant, this definition coincides with the Euler products defined in
\eqref{eqn:EPs}. By general conjectures,
the global $L$-function, which is convergent in some right half-plane, should have
meromorphic continuation to all of $\C$ and satisfy the functional equation\footnote{The
local $\varepsilon$-factors also depend on the choice of a local additive character. We
are assuming a standard choice and hence do not reflect it in the notation. The
global $\varepsilon$-factor is independent of the choice of global additive character.}
\begin{equation}\label{generalfctleq}
 \Lambda(1-s,F,\rho)=\varepsilon(s,F,\rho)\Lambda(s,F,\rho),\qquad\text{where}\quad
 \varepsilon(s,F,\rho)=\prod_{p\leq\infty}\varepsilon_p(s,F,\rho).
\end{equation}
Note that in our case $\varepsilon(s,F,\rho)=\varepsilon_\infty(s,F,\rho)$.
Hence, the functional equations
\eqref{eqn:FEs} are all special cases
of the general conjectured functional equation (\ref{generalfctleq}).
\section{Checking the functional equation}\label{sec:checkFE}

As mentioned in the introduction, the degree 10 adjoint
L-function associated to a Siegel modular form has not
been proven to satisfy a functional equation.
We develop a method of checking a conjectured functional
equation, and in Theorem~\ref{thm:testFE}
we provide a quantitative result that estimates the
probability that this test could yield a false positive.

The main idea behind our method of testing a functional
equation is that an L-function can be evaluated, at a given point
and to a particular accuracy, using finitely many of its
Dirichlet series coefficients.  That evaluation makes fundamental
use of the functional equation.
Furthermore, this can be
done in more than one way.  The consistency of those
calculations provides a check on the  functional
equation.  We quantify the ``probability'' that the
calculations are accidentally consistent by viewing the
coefficients of the $L$-function as a random variable.

In the next section we describe the approximate functional equation
and use it to evaluate an L-function.
Then in Section~\ref{ssec:checkFE} we elaborate on the
ideas from \cite{FKL} to develop our method to check the functional
equation  for the degree-10 Euler product associated to a
Siegel modular form.

\subsection{Smoothed approximate functional equations}\label{sec:appFE}
The material in this section is taken directly from Section~3.2 of~\cite{Rub}.

Let
\begin{equation}
   L(s) = \sum_{n=1}^{\infty} \frac{b_n}{n^s}
\end{equation}
be a Dirichlet series that converges absolutely in a half plane, $\Re(s) > \sigma_1$.

Let
\begin{equation}
    \label{eq:lambda}
    \Lambda(s) = Q^s
                 \left( \prod_{j=1}^a \Gamma(\kappa_j s + \lambda_j) \right)
                 L(s),
\end{equation}
with $Q,\kappa_j \in {\mathbb{R}}^+$, $\Re(\lambda_j) \geq 0$,
and assume that:
\begin{enumerate}
    \item  $\Lambda(s)$ has a meromorphic continuation to all of ${\mathbb{C}}$ with
           simple poles at $s_1,\ldots, s_\ell$ and corresponding
           residues $r_1,\ldots, r_\ell$.
    \item $\Lambda(s) = \varepsilon \cj{\Lambda(1-\cj{s})}$ for some
          $\varepsilon \in {\mathbb{C}}$, $|\varepsilon|=1$.
    \item For any $\sigma_2 \leq \sigma_3$, $L(\sigma +i t) = O(\exp{t^A})$ for some $A>0$,
          as $\abs{t} \to \infty$, $\sigma_2 \leq \sigma \leq \sigma_3$, with $A$ and the constant in
          the `Oh' notation depending on $\sigma_2$ and $\sigma_3$. \label{page:condition 3}
\end{enumerate}

To obtain a smoothed approximate functional equation with desirable
properties, Rubinstein \cite{Rub} introduces an auxiliary function.
Let $g: \C \to \C$ be an entire function that, for fixed $s$, satisfies
\begin{equation*}
    \abs{\Lambda(z+s) g(z+s) z^{-1}} \to 0
\end{equation*}
as $\abs{\Im{z}} \to \infty$, in vertical strips,
$-x_0 \leq \Re{z} \leq x_0$. The smoothed approximate functional
equation has the following form.
\begin{theorem}
    \label{thm:formula}
    For $s \notin \cbr{s_1,\ldots, s_\ell}$, and $L(s)$, $g(s)$ as above,
    \begin{align}
         \label{eq:formula}
         \Lambda(s)  =
g(s)^{-1} \biggl(
         \sum_{k=1}^{\ell} \frac{r_k g(s_k)}{s-s_k}
         + Q^s &\sum_{n=1}^{\infty} \frac{b_n}{n^s} f_1(s,n) \notag \\
         + \varepsilon Q^{1-s} &\sum_{n=1}^{\infty} \frac{\cj{b_n}}{n^{1-s}} f_2(1-s,n)
\biggr)
    \end{align}
    where
    \begin{align}
        \label{eq:mellin}
        f_1(s,n) &:= \frac{1}{2\pi i}
                    \int_{\nu - i \infty}^{\nu + i \infty}
                    \prod_{j=1}^a \Gamma(\kappa_j (z+s) + \lambda_j)
                    z^{-1}
                    g(s+z)
                    (Q/n)^z
                    dz \notag \\
        f_2(1-s,n) &:= \frac{1}{2\pi i}
                    \int_{\nu - i \infty}^{\nu + i \infty}
                    \prod_{j=1}^a \Gamma(\kappa_j (z+1-s) + \cj{\lambda_j})
                    z^{-1}
                    g(s-z)
                    (Q/n)^z
                    dz
    \end{align}
    with $\nu > \max \cbr{0,-\Re(\lambda_1/\kappa_1+s),\ldots,-\Re(\lambda_a/\kappa_a+s)}$.
\end{theorem}

We assume $L(s)$ continues to an entire function, so the first
sum in \eqref{eq:formula} does not appear.  For fixed
$Q,\kappa,\lambda,\varepsilon$, and sequence~$b_n$, and $g(s)$ as described
in the next section, the right side of \eqref{eq:formula}
can be evaluated to high precision. 

We illustrate the approximate functional equation with an example.
Our examples use the genus 2 Siegel modular form which is the unique
weight 20 eigenform, $F$, that is not a Saito-Kurokawa lift.
We consider the  degree 10 adjoint $L$-function associated
to $F$, which we denote by $L_{F,10}$.  Conjecturally
it satisfies the functional equation 
given in \eqref{eqn:FEs} in Proposition~\ref{prop:45101416}.

It is convenient to instead evaluate the ``Hardy function''
$Z_{F,10}$ associated to $L_{F,10}$.  This function is defined by property that $Z_{F,10}(\tfrac12+i t)$ is real
if $t$ is real, and $|Z_{F,10}(\tfrac12+i t)|=|L_{F,10}(\tfrac12+i t)|$.  We use $L$ and $Z$ interchangeably in our discussions. 

If we let $g(s)=1$ and $s=1/2+i$ then \eqref{eq:formula} gives
\begin{align}\label{eqn:ex1}
Z_{F,10}(\tfrac12+i)=\mathstrut & 1.15426 + 0.778012 \,b_2 + 0.50246 \,b_3 + 0.33776\, b_4
+ 0.235813 \,b_5 \cr
&+ \cdots + 0.0000142432 \,b_{82} + 0.0000132692 \,b_{83} + 
\cdots \cr
& +2.8771 \times 10^{-7} \,b_{149} + 2.7402\times 10^{-7} \,b_{150} + \cdots .
\end{align}
Here and throughout this paper, decimal values are truncations
of the true values.
The numerical calculations were done in Mathematica~7.0.1.0
on a Dell Inspiron 9300 laptop running RedHat Linux.

If instead we let $g(s)=e^{-3i s/2}$ and keep $s=1/2+i$ then \eqref{eq:formula}
gives
\begin{align}\label{eqn:ex2}
Z_{F,10}(\tfrac12+i)=\mathstrut &
1.3044 + 0.678149 \, b_{2} + 0.314111 \, b_{3} + 0.12853 \, b_{4} + 
 0.0341584 \, b_{5}
\cr
&+\cdots +0.0000147237 \, b_{82} + 0.0000123925 \, b_{83} +
\cdots \cr
& -1.28515\times 10^{-6} \, b_{149} - 1.22359\times 10^{-6} \, b_{150} 
+ \cdots .
\end{align}

To obtain a numerical value for $Z_{F,10}(\tfrac12+i)$ we need to know
the coefficients $b_n$.
The calculations in~\cite{KohnenKuss} provide
the Satake parameters of $F$ for the primes $p\le 79$,
so Proposition~\ref{prop:45101416} gives the local factors in the
Euler product for~$p\le 79$.  Expanding the product gives
values for infinitely many $b_n$, including all $n\le 82$,
all composite $83\le n \le 79^2$, etc.  
Using the known $b_n$ gives an approximation to $Z_{F,10}(\tfrac12+i)$,
and the numbers in \eqref{eqn:ex1} or \eqref{eqn:ex2} make
it seem plausible that the contribution of the unknown~$b_n$
is small, but we wish to make this precise.  This has two
ingredients:  a bound on $b_n$ and a bound on the terms which
appear in the approximate functional equation.  For the
unknown $b_n$ we assume the Ramanujan bound, so for example
$|b_p|<10$ if $p$ is prime.  For the contribution from the
terms in the approximate functional equation, we use the bound in
Lemma~\ref{lem:booker} for large~$n$, and directly calculate
for smaller~$n$.  
See Section~\ref{sec:booker} for more details.

The results for \eqref{eqn:ex1} and \eqref{eqn:ex2}
are, respectively
\begin{align}\label{eqn:twoZ}
Z_{F,10}(\tfrac12+i)=\mathstrut &3.084662 \pm 0.00047 \cr
Z_{F,10}(\tfrac12+i)=\mathstrut &3.084649 \pm 0.00056 .
\end{align}
The values in \eqref{eqn:twoZ} are consistent with each other.
We view this as a confirmation of the conjectured functional
equation for $L_{F,10}$. 

We summarize the results of similar calculations, for various $s$ and functions
$g$, in Table~\ref{tab:testFE}.  Each column of the table corresponds
to a value for $s$,
and each row corresponds to a function $g(s)=e^{-i \beta s}$
in Theorem~\ref{thm:formula}.
Scanning down each column shows that the values are consistent,
which gives a
check on the functional equation for $L_{F,10}$.
We make this more precise in the next section.

\begin{table}
\begin{tabular}{|c||r@{}l|r@{}l|r@{}l|r@{}l|r@{}l|}
\hline
{$\beta$} &
{$\tfrac12$} & &
$\tfrac12$&+$i$ &
$\tfrac12$&+$2 i$ &
$\tfrac12$&+$3i$ &
$\tfrac12$&+$4i$
\\
\hline\hline
$ 0 $ &2&.148764 &3&.084662 & 3&.263120 & -0&.403124 & 0&.446949   \\
  &\ \ \ \ $ \pm$ 0&.00016 &\ \ \ \ $ \pm$ 0&.00046
 &\ \ \ \   $\pm$ 0&.0044  & \ \ \ \ $ \pm$ 0&.071 & \ \ \ \ $
\pm$ 1&.48   \\
\hline
$\frac14$ &2&.148757 &3&.084643 & 3&.262960 & -0&.405569 &
0&.396311   \\ 
  &\ \ \ \ $ \pm$ 0&.000027 &\ \ \ \ $ \pm$ 0&.00011
  &\ \ \ \   $\pm$ 0&.0013  & \ \ \ \ $ \pm$ 0&.023 & \ \ \ \ $
\pm$ 0&.50   \\  
\hline
$\frac12$ &2&.148743 &3&.084617 & 3&.262768 & -0&.407940 & 0&.356202   \\ 
 &\ \ \ \ $ \pm$ 0&.00021 &\ \ \ \ $ \pm$ 0&.00034     &\ \ \ \   $\pm$ 0&.0019
& \ \ \ \ $ \pm$ 0&.018 & \ \ \ \ $ \pm$ 0&.21  \\  \hline
$ \frac34 $ & 2&.148744 & 3&.084617 & 3&.262767 & -0&.407989 & 0&.355043 \\
&\ \ \ \ $ \pm$ 0&.00014 & \ \ \ \ $ \pm$ 0&.00025
   &\ \ \ \   $\pm$ 0&.014  & \ \ \ \ $ \pm$ 0&.0019  & \ \ \
\ $ \pm$ 0&.16   \\  
\hline
$1$ & 2&.148772 & 3&.0846503 & 3&.262906 & -0&.406974 & 0&.365212   \\   
&\ \ \ \ $ \pm$ 0&.00039 & \ \ \ \ $ \pm$ 0&.00037    &\ \ \ \   $\pm$ 0&.0011
& \ \ \ \ $ \pm$ 0&.0056  & \ \ \ \ $ \pm$ 0&.039   \\
\hline 
$2$ & 2&.148146 & 3&.084355 & 3&.262411 & -0&.408305 & 0&.361331 \\ 
 &\ \ \ \ $ \pm$ 0&.0087 & \ \ \ \ $ \pm$ 0&.0040  
 &\ \ \ \   $\pm$ 0&.0069  &\ \ \ \ $ \pm$ 0&.019  & \ \ \ \ $ \pm$ 0&.071   \\ 
\hline
$3$ & 2&.296591 & 3&.108788 & 3&.277819 & -0&.391079 & 0&.388505   \\   
&\ \ \ \ $ \pm$ 2&.55 & \ \ \ \ $ \pm$ 0&.42   
&\ \ \ \   $\pm$ 0&.25  &\ \ \ \ $ \pm$ 0&.26  & \ \ \ \ $ \pm$ 0&.38   \\
\hline 
\end{tabular} \caption{\sf
\label{tab:testFE} Values obtained for
$Z_{F,10}(\frac12+iT,\beta)$
using equation~\eqref{eq:rearranged} with test function $g(s)=e^{-i \beta s}$,
and using the known Satake parameters
for $F$ for $p\le 79$. }

\end{table}

\subsection{Numerically checking the functional equation}\label{ssec:checkFE}

We wish to check that the
adjoint L-function given by the Euler product~\eqref{eqn:EPs} 
has an analytic
continuation which satisfies the functional equation.
This requires that we evaluate the function outside the region
where the Euler product converges, and also check that 
these values are consistent with the functional equation. 
In what follows we fix the genus 2 Siegel modular form to be the unique
weight 20 eigenform, $F$, that is not a Saito-Kurokawa lift.
Recall that we can use the calculations in~\cite{KohnenKuss} to determine
the Satake parameters of $F$ for the primes $p\le 79$.

Let $g(s)=g(s,\beta) = e^{-i \beta s}$ in Theorem~\ref{thm:formula}.
This meets the conditions of the theorem if
$|\beta|<\frac{\pi}{4}\sum \kappa_j$.  We use 
equation~\eqref{eq:formula} to test the functional equation.
This cannot be done in a naive way, because 
$\Lambda(s)$, as given by the right side of~\eqref{eq:formula},
automatically satisfies $\Lambda(s)=\varepsilon \overline{\Lambda}(1-s)$.
Instead we exploit the fact that the right side~\eqref{eq:formula},
with our choice of $g$, has $\beta$ as a free parameter.

Rewrite \eqref{eq:formula} as
$\Lambda(s) = g(s)^{-1} \Upsilon(s,Q,\kappa,\lambda,\varepsilon,\{b_n\},\beta)$,
and let
\begin{equation}\label{eq:rearranged}
L(s,\beta) = Q^{-s}  \left( \prod_{j=1}^a \Gamma(\kappa_j s + \lambda_j) \right)^{-1}
g(s,\beta)^{-1} \Upsilon(s,Q,\kappa,\lambda,\varepsilon,\{b_n\},\beta).
\end{equation}
Our test for the functional equation of $L(s)$ is that
$L(s,\beta)$ is independent of~$\beta$.
That is, we check the \emph{consistency equation}
\begin{equation}\label{eqn:consistency}
Z(s,\beta_1)-Z(s,\beta_2) = 0.
\end{equation}

For example, using \eqref{eqn:ex1} and \eqref{eqn:ex2} gives
\eqref{eqn:consistency} in the form
\begin{align}
Z_{F,10}(\tfrac12+i,0)-\mathstrut  Z_{F,10}&(\tfrac12+i,3/2)  \cr
= \mathstrut & 0.150138 - 0.0998628 \, b_{2} - 0.188349 \, b_{3} - 0.20923 \, b_{4} - 
 0.201655 \, b_{5} \cr
&+\cdots +4.80503\times10^{-7} \, b_{82} - 8.76677\times10^{-7} \, b_{83} + \cdots \cr
& -1.57286\times10^{-6} \, b_{149} - 1.49761\times10^{-6} \, b_{150} +\cdots \cr
=\mathstrut & 0.
\end{align}

As described immediately before equation~\eqref{eqn:twoZ}, we can estimate
the contribution of the $b_n$ which are not known.  The result
is 
\begin{align}\label{eqn:ex3}
Z_{F,10}(\tfrac12+i,0)-\mathstrut Z_{F,10}&(\tfrac12+i,3/2) \cr
 =\mathstrut & 0.150138 -  0.0998628 \, b_{2} - 0.188349 \, b_{3} - 0.20923 \, b_{4} - 
 0.201655 \, b_{5} \cr
&+\cdots +4.80503\times10^{-7} \, b_{82} +\cdots  - 1.49761\times10^{-6} \, b_{150} +\cdots \cr
 =\mathstrut & \Theta \times 0.00077,
\end{align}
where $|\Theta|\le 1$.

Now we can explain our method of evaluating the strength of
 \eqref{eqn:ex3} as a test of the conjectured functional equation.
We wish to quantify the intuitive
notion that it is unlikely for \eqref{eqn:ex3} to be true just
by chance, because the coefficients of the $b_n$ are large
compared to the
right side of the equation.
We do this by considering the $b_n$ to be random variables,
and furthermore we make some assumptions about their probability
density functions.  This, of course, requires some justification which
we now provide.

L-functions naturally fall into families~\cite{KS,CFKRS},
and the collection of L-functions in a family can be modeled statistically. 
For example, for the family of GL(2) L-functions, each coefficient
$b_j$ has a particular distribution.  The distribution of
$b_p$, for $p$ prime, tends to the Sato-Tate distribution as
$p\to\infty$, and $b_n$, $b_m$ are uncorrelated if $(n,m)=1$.
See~\cite{Ser,CDF} for details.

For other families, there are other distributions,
see~\cite{KeSu} for several examples.  These distributions are
the distributions of traces of matrices in a compact group,
weighted according to Haar measure.
For the Siegel modular forms we
consider here, the Hecke eigenvalues are expected to be distributed
according to an $Sp(4,\Z)$ analogue of the $GL(2)$ case.
This leads to a conjecture for the distribution of the Dirichlet
series coefficients of the degree-10 L-function we are considering
here.  (See additional comments at the end of this section.)

Thus, over the family of L-functions  associated to Siegel modular
forms, we assume the $b_p$ behave as independent random variables, each of
which has a continuous probability distribution which is
supported on $[-10,10]$ and which is bounded
by~$1$, say.  If we focus on one coefficient, say $b_3$, and first
choose all the other $b_n$, then~\eqref{eqn:ex3} becomes
\begin{equation}
C - 0.188349 \, b_{3} = \Theta \times 0.00077,
\end{equation}
where $C$ is some number.
Hence, there is a $C'$ so that,
\begin{equation}\label{eqn:b3prob}
 b_{3} \in\mathstrut  [C'-0.004088,C'+0.004088] .
\end{equation}
Since the PDF of $b_3$ is assumed to be bounded by~$1$,
the probability of \eqref{eqn:b3prob}  being true
is less than the length of that
interval, which is $0.00817$.
In other words, there is less than a 1~percent chance that
$L_{F,10}$ would accidentally pass that test for
satisfying the functional equation.  We have proven:

\begin{theorem}\label{thm:testFE}
Fix the parameters in the functional equation of $L(s)$, as described
at the beginning of Section~\ref{sec:appFE}, and suppose coefficients
$b_j$, $j\in J$ are known and the remaining coefficients obey the
Ramanujan bound.
Let $L(s,\beta)$ be given by \eqref{eq:rearranged},
and choose real numbers $\beta_1$, $\beta_2$ and a complex
number~$s_0$.  Write
\begin{align}\label{eqn:thmtest}
Z(s_0,\beta_1)-Z(s_0,\beta_2) = \mathstrut & \sum_j v_j b_j \cr
= \mathstrut & \sum_{j\in J} v_j b_j  + \Theta \delta ,
\end{align}
where $|\Theta|<1$ and $\delta$ is determined as described
in Section~\ref{sec:booker}.

If the $b_j$ for $j\in J$ are chosen independently from continuous probability
distributions whose PDFs~are bounded by~1, then the probability
that \eqref{eqn:thmtest} is consistent with the functional equation is
less than~$\delta/|v_j|$ for any $j\in J$.
\end{theorem}

It is easy to extend Theorem~\ref{thm:testFE} to the case of several equations.  Suppose we know $b_j$ for $j\in J$, and choose
$s_k$, $\beta_{k,1}$, and $\beta_{k,2}$ for $1\le k\le K$.  We have
\begin{align}\label{eqn:system1}
Z(s_k,\beta_{k,1}) - Z(s_k,\beta_{k,2}) =\mathstrut & \sum_j v_{k,j}b_j \cr
=\mathstrut & \sum_{j\in J}  v_{k,j}b_j + \Theta_k \delta_k .
\end{align}
Select $b_{j_1},\ldots,b_{j_K}$, and suppose all the other
$b_j$ have been determined.  
Then 
the system
\begin{equation}
\{ Z(s_k,\beta_{k,1}) - Z(s_k,\beta_{k,2}) =0\}_{k=1}^K
\end{equation}
is equivalent to
\begin{equation}
A (b_{j_1},\ldots,b_{j_K}) \in (C_1,\ldots,C_K)+[-\delta_1,\delta_1]\times\cdots
		\times [-\delta_K,\delta_K],
\end{equation}
where $A$ is the matrix $(v_{k,j_k})$.  So we can rewrite the condition on
the $b_{j_k}$ as
\begin{equation}\label{eqn:vectorprob}
(b_{j_1},\ldots,b_{j_K}) \in  A^{-1} (C_1,\ldots,C_K)+  A^{-1}([-\delta_1,\delta_1]\times\cdots
                \times [-\delta_K,\delta_K]).
\end{equation}
Since the PDFs of $b_{n_j}$ are assumed to be bounded by~1,
the probability that \eqref{eqn:vectorprob} occurs is bounded
by the volume of the right side, which is
$2^K |\det A|^{-1} \delta_1\cdots \delta_K$.

Here is an example using some of the data from
Table~\ref{tab:testFE}.  
Pairing
the 1st and 5th entries of column~2, and the
5th and 6th entries in column~3, we find:
\begin{align}\label{eqn:samplesystem}
Z(\tfrac12+i,0) - Z(\tfrac12+i,1) =\mathstrut &
-0.07393 + 0.05869 b_2 + 0.10175 b_3 +\cdots \pm 0.00013\cr
Z(\tfrac12+ 2i,1) - Z(\tfrac12+2i,2) =\mathstrut &
-0.41376 + 0.18021 b_2 + 0.43401 b_3 +\cdots \pm 0.0077.
\end{align}

Using the coefficients of $b_2$ and $b_3$ in \eqref{eqn:vectorprob}
we find that the
probability of the above system being satisfied for random
$b_2$, $b_3$ is less than $0.00059$.
This strikes us as rather convincing evidence that the
expected functional equation of $L_{F,10}$ is in fact
correct.

We have a few comments on these calculations.  Our purpose
is to show that it is possible to quantify the precision to which
changes to the test function in the approximate functional equation
give a check on the functional equation of an L-function.
Since it is, in a sense, nonsensical to treat the known coefficients
of an L-function as random, we have tried not to push the analogy
too far.  If the coefficients $b_j$ really were independent and
random, then the sum involving $b_j$ in equations like
\eqref{eqn:ex3} would have a very large variance and the
probability that \eqref{eqn:ex3} holds would be much smaller
than our estimate.  We chose to focus on just one or two
coefficients at a time in order to not stretch plausibility
too much.
Also, the method of equation~\eqref{eqn:vectorprob}
gives poor results if the matrix is close to singular.
In fact, one can add a new equation and \emph{increase}
the probability that the system is consistent, which is
absurd.  This can happen in practice:  adding another equation
based on Table~\ref{tab:testFE} to \eqref{eqn:samplesystem}
actually gives worse results.  This is due to a dependence
among the equations, arising from the fact that for small $t$
it takes relatively few coefficients to evaluate~$L(\tfrac12+i t)$.
We will return to this topic in a subsequent paper.

\subsection{Rigorously evaluating L-functions}\label{sec:booker}

In Section~\ref{sec:appFE} we estimated the contribution of the terms
involving the $b_n$ which were not known explicitly, but only
assumed to satisfy the Ramanujan bound.  This involves estimating
the contribution of infinitely many terms and occurs in two
steps.  First, using Lemma~\ref{lem:booker} below, we determine
$N$ and $\delta_1$ so that the terms involving $b_n$ with $n>N$
contribute, in total, less than $\delta_1$.  Then we explicitly
evaluate the contributions of the terms $f_1(s,n)$ and $f_2(s,n)$
occurring in~\eqref{eq:formula} for $83\le n\le N$;
call that contribution~$\delta_2$.  Then our estimate
for the contribution of the unknown terms is $\delta_1+\delta_2$.

For example, let $\beta=1$ and $s=\frac12 + i$, in order to obtain
the entry in the 5th row and second column of Table~\ref{tab:testFE}.
With $N=10,000$ we find $\delta_1<10^{-6}$ and
$\delta_2 <0.000373$, as reported.  This approach was used to determine
the values in Table~\ref{tab:testFE} and elsewhere
in Section~\ref{sec:checkFE}.  

The following is a very slight modification of
Lemma~5.2 of Booker~\cite{Booker}.

\begin{lemma}\label{lem:booker}
Let
\begin{equation}
G^*(u;\eta,\{\mu_j\}):= \frac{1}{2\pi i} 
\int_{\nu} e^{(u+i\frac{\pi\tau}{4}\eta)(\frac12 - s)}
\prod_{j=1}^r \Gamma_\R(s+\mu_j) \frac{ds}{s}.
\end{equation}
Then for $X\ge r$,
\begin{equation}
G^*(u;\eta,\{\mu_j\}) \le \frac{K r}{X} e^{\Re(\mu) u} e^{-X}
\prod_{j=1}^r \left(1+\frac{r\nu_j}{X} \right)^{\nu_j} ,
\end{equation}
where
$\delta=\frac{\pi}{2}(1-|\eta|)$, $\nu_v=\frac12(\Re\mu_j-1)$,
$\mu=\frac12+\frac{1}{r}(1+\sum\mu_j)$,
$K=2\sqrt{\frac{2^{r+1}}{r}\frac{e^{\delta(r-1)}}{\delta}}
e^{-\frac{\pi r \eta \Im\mu}{4}}$,
and $X=\pi r \delta e^{-\delta}e^{2u/r}$.
\end{lemma}

Note that our $G^*$ is identical to the function $G$ in
 Lemma~5.2~of~\cite{Booker} except for the extra factor
of $1/s$ in the integrand.

\begin{proof}
Move the line of integration to the $2\sigma$ line
and let $s=2\sigma +2 i t$ and use the
trivial estimate $1/|s| \le 1/(2\sigma)$ to get
\begin{equation}
G^*(u;\eta,\{\mu_j\}) \le \frac{1}{\sigma}
\frac{1}{2\pi} 
\int_{2\sigma} \left| e^{(u+i\frac{\pi\tau}{4}\eta)(\frac12 - s)}
\prod_{j=1}^r \Gamma_\R(s+\mu_j)\right| {dt}.
\end{equation}
Now exactly follow the proof of Lemma~5.2~in~\cite{Booker},
which in the last step chooses $\sigma=X/r$.
\end{proof}

\bibliographystyle{plain}
\bibliography{degree_ten_pac_j_final}

\end{document}